\documentclass[12pt]{article}
\usepackage{graphicx}
\usepackage{epsfig}
\usepackage{amsmath}
\usepackage{amssymb}
\usepackage{amsthm}
\usepackage{latexsym}
\usepackage{cite}
\usepackage{amsbsy,amsfonts,euscript,eepic,rotating}

\newcommand{\comment}[1]{}

\newtheorem{lemma}{Lemma}[section]

\newtheorem{corollary}{Corollary}[section]
\newtheorem{proposition}{Proposition}[section]

\newtheorem{hypothesis}{Hypothesis}[section]

\begin{document}

\title{\LARGE
%{\bf Cluster Area Measures and the \\ Ising (Conformal) Field Theory}
%{\bf The Ising (Conformal) Field Theory \\ as a Random Signed Measure}
%{\bf Ising (Conformal) Fields \\ and Random Signed Measures}
{\bf Ising (Conformal) Fields \\ and Cluster Area Measures}
}

\author{
{\bf Federico Camia}
\thanks{Research supported in part by a Veni grant of the
NWO (Dutch
Organization for Scientific Research).}\,
%and by the U.S. NSF under
%grant PHY99-07949.}\,
\thanks{E-mail: fede@few.vu.nl}\\
{\small \sl Department of Mathematics, Vrije Universiteit Amsterdam}\\
{\small \sl De Boelelaan 1081a, 1081 HV Amsterdam, The Netherlands}
\and
{\bf Charles M.~Newman}
\thanks{Research supported in part by the NSF under grants
DMS-06-06696 and OISE-07-30136.}\,
\thanks{E-mail: newman@courant.nyu.edu}\\
{\small \sl Courant Inst.~of Mathematical Sciences,
New York University}
\\
{\small \sl 251 Mercer Street, New York, NY 10012, USA}
%\\
%{\small \sl http://math.nyu.edu/faculty/newman/index.html}
%, New York, NY 10012, USA}
}

\date{}

\maketitle

\begin{abstract}
We provide a representation for the scaling limit of the $d=2$
critical Ising magnetization field as a (conformal) random field
using SLE (Schramm-Loewner Evolution) clusters and associated
renormalized area measures. The renormalized areas are from the
scaling limit of the critical FK (Fortuin-Kasteleyn) clusters and
the random field is a convergent sum of the area measures with
random signs. Extensions to off-critical scaling limits, to $d=3$
and to Potts models are also considered.
\end{abstract}



\section{Introduction} \label{sec-intro}

The Ising model in $d=2$ dimensions is perhaps the most studied
statistical mechanical model and has a special place in the theory
of critical phenomena since the groundbreaking work of Onsager~\cite{onsager}.
Its scaling limit at or near the critical point is recognized to
give rise to Euclidean (quantum) field theories. In particular,
the scaling limit of the lattice magnetization field should be a
Euclidean random field and, at the critical point, the simplest
reflection-positive conformal field theory
$\Phi^0$~\cite{bpz2,cardy08}. As such, there have been a variety
of representations in terms of free fermion fields~\cite{sml} and
explicit formulas for correlation functions (see, e.g.,
\cite{mw-book,palmer} and references therein). In this paper, we provide
a construction of $\Phi^0$ in terms of random geometric
objects associated with Schramm-Loewner Evolutions (SLEs)~\cite{schramm}
(see also~\cite{cardy3,kn,lawler2,werner4}) and Conformal Loop Ensembles
(CLEs)~\cite{sheffield,shw,werner3} --- namely, a gas (or random
process) of continuum loops and associated clusters and (renormalized)
area measures.

Two such loop processes arise in the results announced by
Smirnov~\cite{smirnov-icm,smirnov-ising1,smirnov-ising2,smirnov-ising3,smirnov-ising4}
(see also the work of Riva and Cardy in~\cite{rc} --- in particular Sections
6 and 7 there) that the full scaling limit of critical Ising spin cluster
boundaries (respectively, FK random cluster boundaries) is given by
the (nested version of) CLE with parameter $\kappa=3$ (resp.,
$\kappa=16/3$).
%We show that
One can try to associate with each continuum cluster ${\cal C}_j^*$ %{[\kappa]}$
or external boundary loop ${\cal L}_j^*$ in  %{[\kappa]}$ in
the scaling limit a finite area measure $\mu_j^*$ representing the
rescaled number of sites in the corresponding lattice cluster (where
$*$ is SP for the spin case and FK for the random cluster case).
%%$\kappa=3$ for the spin case and $16/3$ for the FK case).
We can in fact do this for the FK case and expect it to also be valid
for the spin case.

Although one might try to represent the Euclidean field $\Phi^0$
using spin clusters
%and CLE$_3$
by a sum $\sum_k \chi_k \mu_k^{SP}$, where %{[3]}$,
the $\chi_k$'s are
%deterministically
$+1$ or $-1$
depending on whether ${\cal C}^{SP}_k$ corresponds to a $+$ or
%{[3]}_k$ corresponds to a $+$ or
$-$ spin cluster, this does not seem to work. Instead, we use the
FK clusters, which leads to $\Phi^0 = \sum_j \eta_j
\mu_j^{FK}$, where the $\eta_j$'s are %{[16/3]}$, where the $\eta_j$'s are
independent random signs.
The (countable) family $\{\mu_j^{FK}\}$ is a ``point'' process with
each $\mu_j^{FK}$ a ``point'' and where distinct ``points'' should
be orthogonal measures.

For a bounded $\Lambda \subset {\mathbb R}^2$
with nonempty interior, one expects that
$\sum_j \mu_j^{FK} (\Lambda) = \infty$.
%The expectation of divergence has two sources.
%One is that this would
This would follow from the scaling covariance
expected for $\{\mu_j^{FK}\}$ and described at the end
of this section.
%A second source is that this is the case
The same happens for the corresponding measures in independent
percolation that count so-called ``one-arm'' sites, as follows from
work of Garban, Pete and Schramm~\cite{gps,garban}. Nevertheless,
%({\tt don't we actually knwon
%this from the fact that without random signs the unnormalized sum grows
%like ${1/a}^2$ while $\Theta_a$ goes like $a^{15/8}$?}),
for any $\varepsilon >0$
only finitely many $\mu_j^{FK}$'s will have support
that intersects  $\Lambda$
{\it and\/} has diameter greater than~$\varepsilon$.
Furthermore, with probability one,
$\sum_j [\mu_j^{FK} (\Lambda)]^2 < \infty$ which leads to convergence
(at least in $L_2$) of the sum with random signs
$\sum_j \eta_j \mu_j^{FK}(\Lambda)$.
%One feature of this representation is that
%the sum is convergent with random signs, but presumably not without.
We note that divergence of $\sum_j \mu_j^{FK}$ means that $\Phi^0 = \sum_j \eta_j
\mu_j^{FK}$ is not a signed measure; i.e., even restricted to a bounded
$\Lambda$, it is not the difference of two positive {\it finite\/} measures.
For negative results of a similar sort, but in the context of Gaussian
random fields, see~\cite{colella-lanford}.

In the next section, we set up notation for the Ising model on the
square lattice and its FK representation and review how the scaling
limit of FK cluster boundaries may be viewed as a
process of noncrossing continuum loops ${\cal L}^{FK}_j$ and %[\kappa]}_j$ and
associated continuum clusters ${\cal C}^{FK}_j$. We then show
why the natural scaling for the Ising spin variables at criticality
to obtain a Euclidean (random) field $\Phi^0$ leads %in the FK case
%($\kappa=16/3$)
to natural rescaled area measures ${\mu}^{FK}_j$ %[\kappa]}_j$
supported on ${\cal C}^{FK}_j$ %[\kappa]}_j$
and to the representation of
$\Phi^0$ in terms of those measures. We also discuss why area measures
${\mu}^{SP}_k$ for spin clusters are not appropriate for representing
$\Phi^0$, by using an example taken from the infinite temperature
Ising model on the triangular lattice, ${\mathbb T\/}$.

In Section~\ref{sec-area}, we use (see Proposition~\ref{lemma-zero-limit})
and then
discuss how to verify a decay property of the critical Ising two-point
correlation or equivalently the FK connectivity function. Another
essential ingredient in our analysis is a bound (see
Prop.~\ref{lemma-crossprob}) on
the number of macroscopic FK
clusters. Although we focus on
%verifying the decay property for
critical Ising-FK percolation on ${\mathbb Z}^2$, similar arguments
can be applied to other lattices and to independent (and,
in principle, Ising spin)
percolation. The case of independent percolation is
discussed at the end of Section~\ref{sec-area}.
%In the context of general translation
%invariant critical percolation models
%on regular lattices, we give sufficient conditions for this
%property to be valid, verifiable
%in critical FK models both for $q=2$ (Ising) and $q=1$
%(independent percolation). Finally,
In Section~\ref{sec-discussion}, we review the general
conclusions of our work and discuss extensions to off-critical
(or as they are sometimes called, near-critical)
scaling limits, either as temperature $T \to T_c$, the critical
temperature, with magnetic field $h=0$, or else as $h \to 0$ with $T=T_c$.
Finally, we propose there that a cluster area measure representation
should also be valid for the $d=3$ Ising model and for the $d=2$ $q$-state
Potts model with $q=3$ or $4$.
%Our focus in this paper is on \emph{representations} of the scaling
%limit rather than its existence. In particular,

Before concluding this section, we wish to emphasize that this paper
is meant to serve as an introduction, readable by both mathematicians
and physicists, to a \emph{representation} for the Ising scaling limit
field $\Phi^0$ in terms of the limit rescaled area measures $\{\mu_j^{FK}\}$.
We hope this will prove useful in providing a general conceptual framework
for field-based scaling limits like Aizenman and Burchard~\cite{ab} did
for connectivity-based ones. Although detailed explanations and proofs are
provided in this paper for certain issues, others are avoided. In particular,
although the next two sections of the paper provide arguments
%needed to guarantee the
for the existence of both $\Phi^0$ and $\{\mu_j^{FK}\}$
as (subsequence) limits
of the corresponding lattice quantities, they
do not provide the
tools needed to prove that the limits are unique.
This will be done in a future paper in collaboration
with C.~Garban, along with a proof of related
properties such as that $\Phi^0$ and $\{\mu_j^{FK}\}$
have the expected conformal covariance including
%the scaling property
that for $\alpha >0$, $\alpha^{1/8} \Phi^0(\alpha z)$ and
$\{\alpha^{-15/8}\mu_j^{FK}(d(\alpha z))\}$ are equidistributed with
$\Phi^0(z)$ and $\{\mu_j^{FK}(dz)\}$.

\section{Ising (Euclidean) Field } \label{sec-ising-field-theory}

We consider the standard Ising model on the square lattice ${\mathbb Z}^2$
with Hamiltonian
\begin{equation} \label{hamiltonian} \nonumber
{\bf H} = -\sum_{\{x,y\}} S_x S_y - h \sum_x S_x,
\end{equation}
where the first sum is over nearest neighbor pairs in ${\mathbb Z}^2$
(or bonds $b=\{x,y\}$), %the second sum is over $x \in {\mathbb Z}^2$
the spin variables $S_x$ are $(\pm 1)$-valued and the external field
$h$ is in $\mathbb R$.

When there is a unique infinite volume Gibbs distribution for some
value of $h$ and inverse temperature $\beta=1/T$, we denote by
$\langle \cdot \rangle_{\beta,h}$ its expectations. There is a
critical $\beta_c$ such that nonuniqueness occurs only for $h=0$ and
$\beta>\beta_c$. In particular, the critical Gibbs measure is unique
and in that case we use the notation $\langle \cdot
\rangle_c=\langle \cdot \rangle_{\beta_c,0}$. By translation
invariance, the two-point correlation $\langle S_x S_y
\rangle_{\beta,h}$ is a function only of $y-x$, which in the
critical case we denote by $\tau_c(y-x)$.

We want to study the random field associated with the spins on the
rescaled lattice $a{\mathbb Z}^2$ in the scaling limit $a \to 0$.
More precisely, for test functions $f(z)$ of bounded support on
${\mathbb R}^2$, we can define for the critical model
\begin{equation} \label{eq-discrete-field}
\Phi^a(f) = \int_{{\mathbb R}^2} f(z) \Phi^a(z) dz =
\int_{{\mathbb R}^2} f(z) [ \Theta_a \sum_{x
\in {\mathbb Z}^2} S_x \delta(z-ax) ] dz = \Theta_a \sum_{z \in
a{\mathbb Z}^2} f(z) S_{z/a},
\end{equation}
with an appropriate choice of the scale factor $\Theta_a$. Since
$\Phi^a(f)$ is a random variable with zero mean, it is natural to
choose $\Theta_a$ so that $\langle [\Phi^a(f)]^2 \rangle_c$ is
bounded away from 0 and~$\infty$ as $a \to 0$. Choosing $\Theta_a$
so that this second moment is exactly one for $f$ the indicator
function of the unit square $[0,1]^2$ yields
\begin{equation} \label{eq1-Theta}
\Theta_a^{-1} = \sqrt{ \sum_{z,w \in \Lambda_{1,a}} \langle S_{z/a}
S_{w/a} \rangle_c } = \sqrt{ \sum_{x,y \in \Lambda_{1/a}} \tau_c(y-x) },
\end{equation}
where $\Lambda_{L,a} = [0,L]^2 \cap a{\mathbb Z}^2$ and
$\Lambda_L=\Lambda_{L,1} = [0,L]^2 \cap {\mathbb Z}^2$.

One way to formulate the FK representation of the Ising model (for
$h=0$ and $\beta \leq \beta_c$) is that coexisting with the $(\pm
1)$-valued spin variables $S_x$ on the sites $x$ of ${\mathbb Z}^2$
are $\{0,1\}$-valued occupation variables $n_b$ on the bonds $b=\{x,y\}$
of ${\mathbb Z}^2$. The occupied or open ($n_b=1$) FK bonds determine FK
clusters, $C_i$, which are the sets of sites $x$ in ${\mathbb Z}^2$
connected to each other by paths of open FK bonds. One can
generate the $S_x$'s from the $n_b$'s by assigning independent
symmetric $\pm 1$ random signs $\eta_i$ to the $C_i$'s and then
setting $S_x=\eta_i$ for every $x$ in $C_i$. If we write $x
\stackrel{FK}{\longleftrightarrow} y$ to denote that $x$ and $y$ are
in the same FK cluster, it is immediate that
the FK connectivity function at criticality is simply given by
\begin{equation} \label{eq-connect-function} \nonumber
P(x \stackrel{FK}{\longleftrightarrow} y) = \langle S_x S_y
\rangle_c = \tau_c(y-x).
\end{equation}
Denoting by $E_c$ expectation in the critical system, by
$\hat{C}^a_i$ the restriction of the cluster $aC_i$ in $a{\mathbb
Z}^2$ to $[0,1]^2$, and by $|\hat{C}^a_i|$ the number of
($a{\mathbb Z}^2$)-sites in
$\hat{C}^a_i$, we have
\begin{equation} \label{eq2-Theta}
\Theta^{-2}_a = \sum_{x,y \in \Lambda_{1/a}} \tau_c(y-x) = \sum_{x,y
\in \Lambda_{1/a}} P (x \stackrel{FK}{\longleftrightarrow} y) =
E_c(\sum_i |\hat{C}^a_i|^2).
\end{equation}

By the definition of $\Theta_a$ we see that the rescaled areas
$W^a_i = \Theta_a |\hat{C}^a_i|$ are uniformly square summable in
the sense that $E_c \sum_{i} (W^a_i)^2 = 1$ for all $a$. We would
like to argue that, at least along subsequences of $a$'s tending to
zero, $\{ W^a_i \}$ has a nontrivial limit in distribution. This is
already partly clear --- i.e., no $W^a_i$ can diverge to $+\infty$.
But what prevents them all from tending to zero as $a \to 0$? It
turns out that this uses the following hypothesis about
$\tau_c(y-x)$ (where $1/\sqrt{2}$ is the appropriate constant for
the lattice ${\mathbb Z}^2$) --- roughly speaking, that it decays
like $||y-x||^{-2 \theta}$ with $\theta < 1$, where $||\cdot||$
denotes Euclidean norm or that $\sum_{||x||\leq r}\tau_c(x)$
diverges as a power when $r \to \infty$. (It also uses that the
crossing probability of an annulus is bounded away from one as $a
\to 0$ --- see~(\ref{eq-crossprob2}) and
Prop.~\ref{lemma-crossprob}.)

\begin{hypothesis} \label{hypo}
%For some fixed $\theta<1$ and any small $\varepsilon >0$, there is
%some $K>0$ and $K'<\infty$ such that for any $x \in {\mathbb Z}^2$
%with large $||x||$ and any$x_{\varepsilon} \in {\mathbb Z}^2$ with
%$||x_{\varepsilon} - \varepsilon x|| \leq 2$,
%
%$\tau_c(x) \geq K
%\varepsilon^{2\theta} \tau_c(x_{\varepsilon})$.
%
For some fixed $\theta<1$, there are constants $K_1>0$ and
$K_2<\infty$
%such that for any $x \in {\mathbb Z}^2$ with large $||x||$ and any
%$\varepsilon$ small,
such that for any small $\varepsilon >0$ and then for any
$x \in {\mathbb Z}^2$ with large $||x||$,
\begin{equation} \label{eq-connectivity}
K_2 \tau_c(x_{\varepsilon}) \, \geq \, \tau_c(x) \, \geq \,
K_1 \, \varepsilon^{2\theta} \tau_c(x_{\varepsilon}) \,
\end{equation}
for any $x_{\varepsilon} \in {\mathbb Z}^2$ with
$||x_{\varepsilon} - \varepsilon x|| \leq 1/ \sqrt{2}$.
\end{hypothesis}

As we will discuss, the clusters $\{ C^a_i = a C_i \}$ on the
rescaled lattice $a{\mathbb Z}^2$ will converge in the scaling limit
to full plane continuum clusters $\{ {\cal C}_j^{FK} \}$ in ${\mathbb
R}^2$. In that limit most of the lattice clusters disappear because
they are not of macroscopic size. The importance of the lower bound
on $\tau_c(x)$ in~(\ref{eq-connectivity})
%Hypothesis~\ref{hypo} (see Prop.~\ref{lemma-zero-limit}) is that it
is that it guarantees (see Prop.~\ref{lemma-zero-limit})
that the rescaled areas of the microscopic clusters are
negligible (at least in a square summable sense). That is, the
contribution to $\sum_i (W_i^a)^2$ coming from clusters $C^a_i$
whose intersection with the unit square has small macroscopic diameter
is negligible. A corresponding statement is true for the clusters
that contribute to the field $\Phi^a(f)$
for
more general test functions $f$ of bounded support.
The significance of the upper bound on $\tau_c(x)$ in~(\ref{eq-connectivity})
is that it easily implies that $\langle [\Phi^a (f)]^2 \rangle_c$
is bounded away from $0$ and~$\infty$ as $a \to 0$.

In a series of papers, the authors constructed~\cite{cn} a
certain process of loops in the plane and proved~\cite{cn1}
(see also~\cite{cn08})
convergence to it in the scaling limit of the collection of
boundaries of all (macroscopic) clusters for critical
independent site percolation on the triangular lattice.
In the limit there is no
self-crossing or crossing of different loops but there is
self-touching and touching between different loops.
Moreover, the loops are locally SLE$_6$ curves.

\begin{figure}[!ht]
\begin{center}
\includegraphics[width=10cm]{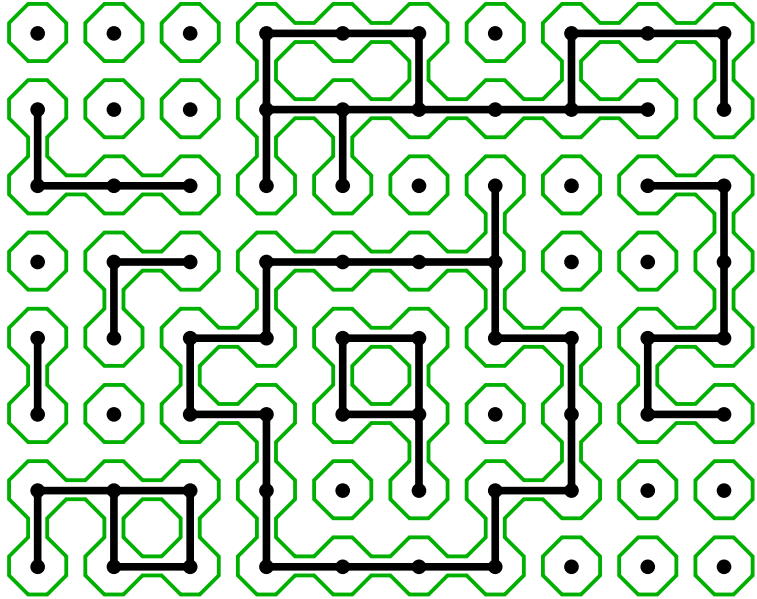}
\caption{Example of an FK bond configuration in a rectangular region and
the associated loops in the medial lattice. Black dots represent sites
of ${\mathbb Z}^2$, black horizontal and vertical edges represent open
FK bonds and the lighter (green) loops are on the medial lattice. We focus
on those loops which have ${\mathbb Z}^2$-sites immediately on their {\it inside\/}.}
\label{loopsfigure}
\end{center}
\end{figure}

Similar results for the 2D critical Ising model
%and (Ising-)FK percolation
on the square lattice have been announced by
Smirnov~\cite{smirnov-icm,smirnov-ising1,smirnov-ising2,smirnov-ising3,smirnov-ising4}.
%In the case of the Ising model
There one considers either the boundaries between plus and minus
spin clusters~\cite{smirnov-ising4}, or the loops in the medial
lattice that separate FK from dual FK clusters~\cite{smirnov-ising3}
(see Figure~\ref{loopsfigure}).
We will focus on those loops which separate FK clusters in
the original ${\mathbb Z}^2$ lattice on their {\it inside\/}
from dual FK clusters in the dual lattice on their {\it outside\/}.
In the scaling limit of spin cluster
boundaries one would obtain simple loops that do not touch each other and
locally are SLE$_3$ type curves. In the case of FK cluster
boundaries, there would instead be
self-touching and touching between
different loops (but no crossing),
like in the percolation case. Now, however, the
loops would locally be SLE$_{16/3}$ type curves.

In the FK case, each loop $L_i^a$ that we consider on the medial lattice
of $a{\mathbb Z}^2$ is the outer boundary of a rescaled FK cluster
$C_i^a$. The inner boundary of $C_i^a$ is made of ``daughter'' loops
$L_{i,k}^a$ corresponding to the ``holes'' in $C_i^a$. In the
scaling limit $a \to 0$, one can analogously identify a continuum
cluster ${\cal C}_j^{FK}$ as the closed set left after removing from
${\mathbb R}^2$ the (open) exterior of the loop ${\cal L}_j^{FK}$
and the (open) interiors of its daughter loops ${\cal L}_{j,n}^{FK}$
(with interiors and exteriors defined using winding numbers). We
remark that because the scaling limit is only a  limit in
distribution and no special effort was made to coordinate indexing
for clusters in the lattice and in the continuum, we use different
letters, $i$ and $j$, for the two indices. We denote by $\{
\mu_j^{FK} \}$ the finite measures supported on $\{ {\cal C}_j^{FK}
\}$ corresponding to the limit of the rescaled areas $\{ W_i^a \}$
as $a \to 0$, in the sense, e.g., that $\{\mu_j^{FK}(\Lambda_1)\}$
is the scaling limit of the rescaled areas $\{W^a_i\}$. The
existence and nontriviality of $\{\mu_j^{FK}(\Lambda_1)\}$ (or
of $\{\mu_j^{FK}(f)= \int f(z) \mu_j^{FK}(dz)\}$ for more general test
functions $f(z)$ of bounded support) will follow from
Hypothesis~\ref{hypo} (see Prop.~\ref{lemma-zero-limit})
and Prop.~\ref{lemma-crossprob}, as noted above.
The collection $\{\mu_j^{FK}\}$ ought to be a functional
of $\{{\cal L}_j^{FK}\}$ %(or of $\{{\cal C}_j^{FK}\}$)
as has recently been proved in the independent percolation context
by Garban, Pete and Schramm~\cite{gps,garban}.

%The $\{({\cal L}_i^0,\{{\cal L}_{i,j}^0 \},\mu_i^0)\}$ process will
%be $\sigma$-finite in that there will be a finite number of $({\cal
%L}_i^0,\{{\cal L}_{i,j}^0 \},\mu_i^0)$'s satisfying two
%restrictions: (i) ${\cal L}_i^0$ touches (or alternately, is
%contained in) a given bounded region of ${\bf R}^2$, and (ii) the
%diameter $diam({\cal L}_i^0)> \varepsilon$. We conclude this
%subsection by noting that the square summability of ${\hat A}_i =
%\mu_i^0(\Lambda_1)$ (or more generally of $\mu_i^0(f)=\int f(z)
%\mu_i^0(dz)$) will be important for the random field construction
%below.

Letting $\{ \eta_j \}$ denote i.i.d.\ symmetric $(\pm 1)$-valued
variables, one obtains the following representation of the Euclidean
field $\Phi^0$: for test functions $f(z)$ of bounded support,
\begin{equation} \label{eq:eucfield}
\Phi^0 (f) \, = \, \sum_j \eta_j \mu_j^{FK} (f) \, = \, \int_{{\mathbb
R}^2} f(z) \, \sum_j \eta_j \mu_j^{FK}(dz) \, .
\end{equation}
To be more precise, the sums in~(\ref{eq:eucfield}) should first be
restricted to clusters with diameter greater than $\varepsilon$
and then convergence (in $L_2$) as the cutoff $\varepsilon \to 0$ will
follow from the square summability discussed earlier.
%where the convergence (in $L_2$) of the sum over $j$ in~(\ref{eq:eucfield})
%%(in $L_2$ of the full probability space for the $\mu_j^{FK}$'s and $\eta_j$'s)
%follows from the square summability discussed earlier.

As noted in the Introduction, one might be tempted to represent the
Euclidean field using spin clusters and hence SLE$_3$ type loops.
%by a summation $\sum_i \chi_i \mu_i^{\kappa}$, where $\kappa=3$ and the
%$\chi_i$'s are deterministically $+1$ or $-1$ depending on whether
%${\cal C}^{\kappa}_i$ corresponds to a $+$ or $-$ spin cluster, this
%does not seem to work.
If we use $\{ \mu_k^{SP+}\}$ and $\{ \mu_{k'}^{SP-} \}$ to denote the
limits of appropriately rescaled areas of plus and minus spin
clusters, respectively, then on a formal level, by decomposing the
righthand side of~(\ref{eq-discrete-field}) into the contribution
from plus and minus clusters, one might expect that $\Phi^0$
of~(\ref{eq:eucfield}) would also be given by $\sum \mu_k^{SP+} - \sum
\mu_{k'}^{SP-}$ (with some resummation needed to handle the difference
of two presumably divergent series) as an alternative to $\sum \eta_j
\mu_j^{FK}$. This appears not to be so, as can be understood by
considering the simple situation of the Ising model on the triangular
lattice~$\mathbb T$ at $\beta=0$.

The latter is noncritical as an Ising model and the correct
Euclidean field obtained by using the noncritical FK clusters (which
are just isolated sites since $\beta=0$) and the
$\beta =0$ version of~(\ref{eq1-Theta}), is two-dimensional Gaussian
white noise. But if one considers the Ising spin clusters, this is
critical independent site percolation on $\mathbb T$ and the formal
expression $\sum \mu_k^{SP+} - \sum \mu_{k'}^{SP-}$, besides the
resummation issue, seems unrelated to white noise.
%while it probably can be made sense of (via a cutoff on
%cluster diameter),
Indeed, applying and then removing a cutoff, as explained
after~(\ref{eq:eucfield}), in this case would probably not lead to
the physically correct limit.

\section{Area Measure} \label{sec-area}

In the previous section we gave a representation of the Ising
Euclidean spin field in terms of rescaled counting measures that
give the ``areas'' of macroscopic Ising-FK clusters. In this
section we first explain how to use Hypothesis~\ref{hypo} to get the
existence of nontrivial limits in distribution of
these area measures, at least along
subsequences of $a$'s tending to zero.
%tightness for the distribution of
%those area measures in the scaling limit $a \to 0$.
%We will work in the context of general translation
%invariant (critical) percolation models on regular lattices (rescaled
%by $a$), including critical FK models. As we will explain later
%in the section,
We then explain how to verify
Hypothesis~\ref{hypo} first for critical Ising-FK percolation
on ${\mathbb Z}^2$ and then for critical independent site or
bond percolation on ${\mathbb T}$ or ${\mathbb Z}^2$.
Using the notation introduced in
Section~\ref{sec-ising-field-theory} and denoting by
$diam(\hat{C}^a_i)$ the Euclidean diameter of $\hat{C}^a_i$, we have
the following proposition.
\begin{proposition} \label{lemma-zero-limit}
Hypothesis~\ref{hypo} implies that
\begin{equation} \nonumber
\lim_{\varepsilon \to 0} \limsup_{a \to 0}
\Theta^2_a E_c(\sum_{i:diam(\hat{C}^a_i)
\leq \varepsilon} |\hat{C}^a_i|^2) = 0.
\end{equation}
\end{proposition}
The usefulness of Prop.~\ref{lemma-zero-limit} is based on the
additional result that for every fixed $\varepsilon$, in the scaling limit there will
only be finitely many FK clusters with diameter larger than $\varepsilon$
that intersect $[0,1]^2$; this important feature of the scaling
limit will be discussed below --- see Prop.~\ref{lemma-crossprob}.
Once one has Prop.~\ref{lemma-crossprob}, it then follows from
Prop.~\ref{lemma-zero-limit} that %$\{ \Theta_a |\hat{C}^a_i|: i=1,2,\ldots \}$
the collection $\{ \Theta_a |\hat{C}^a_i|\}$
has nontrivial subsequential limits; 
i.e., it is not possible that all
$\Theta_a |\hat{C}^a_i|$'s scale to zero as $a \to 0$.
Said more physically, Prop.~\ref{lemma-zero-limit}
implies that there is a negligible contribution to the magnetization
from FK clusters whose linear size is small on a macroscopic lengthscale 
while Prop.~\ref{lemma-crossprob} says that there are only finitely
many larger clusters touching any bounded region. 
Together, they lead to the representation~(\ref{eq:eucfield}). \\

\noindent{\bf Proof of Proposition~\ref{lemma-zero-limit}.}
%Using~(\ref{eq2-Theta}) and Hypothesis~\ref{hypo}, we have
Using Hypothesis~\ref{hypo}, we can compare $\sum_{z' \in
\Lambda_{\varepsilon' r}} \tau_c(z')$ for small $\varepsilon'$ as
$r\to \infty$ to $\sum_{z \in \Lambda_r}\tau_c(z)$ by using the
second inequality of~(\ref{eq-connectivity}) to compare each
$\tau_c(z')$ to the $\tau_c(z)$'s with $\varepsilon' z$ in the unit
length square centered on $z'$ (so that we may take $z'$ as
$z_{\varepsilon'}$). Since there are approximately
$(1/{\varepsilon'})^2$ such $z$ sites, we have that
\begin{equation} \label{eq-taucompare} \nonumber
\liminf_{r \to \infty}
%\{\sum_{z \in \Lambda_r}\tau_c(z)\} /
%\{(1/{\varepsilon'})^2 K_1 ({\varepsilon'})^{2\theta}
%\sum_{z' \in \Lambda_{\varepsilon' r}} \tau_c(z') \} \, \geq \, 1 \, .
\frac{\sum_{z \in \Lambda_r}\tau_c(z)}
{(1/{\varepsilon'})^2 ({\varepsilon'})^{2\theta}
\sum_{z' \in \Lambda_{\varepsilon' r}} \tau_c(z') }
\, \geq \, K_1 \, .
\end{equation}
Using this lower bound (with $r = 1/2a$ and $\varepsilon' = 2 \varepsilon$)
and~(\ref{eq2-Theta}), we have that
\begin{eqnarray}
\limsup_{a\to0} \Theta_a^{2} E_c(\sum_{i:diam(\hat{C}^a_i) \leq \varepsilon}
|\hat{C}^a_i|^2) & \leq & \limsup_{a\to0}
\frac{\sum_{x,y \in \Lambda_{1/a}, ||x-y|| \leq
\varepsilon/a} \tau_c(y-x)}{\sum_{x,y \in \Lambda_{1/a}}
\tau_c(y-x)} \nonumber \\
& \leq & \limsup_{a\to0} \frac{K' (1/a)^2
%\sum_{||z|| \leq \varepsilon/a}
\sum_{z' \in \Lambda_{\varepsilon /a}}
\tau_c(z')}{K'' (1/a)^2
%\tau_c(z)}{K'' (1/a)^2
%\sum_{||z|| \leq 1/2a} \tau_c(z)} \nonumber \\
\sum_{z \in \Lambda_{1/2a}} \tau_c(z)} \nonumber \\
%& \leq & \limsup_{a\to0} \frac{K'
%%\sum_{||z|| \leq \varepsilon/a} \tau_c(z)}{K'''
%\sum_{z' \in \Lambda_{\varepsilon /a}} \tau_c(z')}{K'''
%(1/ \varepsilon)^2 \varepsilon^{2\theta}
%%\sum_{||z|| \leq \varepsilon/a} \tau_c(z)} \nonumber \\
%\sum_{z' \in \Lambda_{\varepsilon /a}} \tau_c(z')} \nonumber \\
& = & K'''' \varepsilon^{2(1-\theta)}. \label{eq-theta-bound}
\end{eqnarray}
The proposition follows from the observation that
the last expression in~(\ref{eq-theta-bound})
tends to zero as $\varepsilon \to 0$ since $\theta<1$. \fbox{} \\

The next two lemmas will be used to verify Hypothesis~\ref{hypo}.
Let $B_x(r)$ denote $\{y \in {\mathbb Z}^2: ||y-x|| \leq r\}$,
%${\mathbb Z}^2 \cup \{z \in {\mathbb R}^2: ||z-x|| \leq r\}$,
and denote its ${\mathbb Z}^2$-boundary
by $\partial B_x(r)$. If the subscript
is omitted, we refer to the disc centered at the origin $0$. We
denote by $P_{\partial B(r)}^W$ ($W$ for wired) the
critical FK measure inside
$B(r)$ with wired (i.e., everything connected) boundary condition on
$\partial B(r)$.
The next lemma is based on the FKG inequalities.
\begin{lemma} \label{lemma-upper-bound}
%For all percolation measures satisfying FGK (e.g., for all FK
%measures with $q \geq 1$ and all $p \in [0,1]$),
\begin{equation} \nonumber
\tau_c(y-x) \leq P_{\partial B(||x-y||/3)}^W (0
\stackrel{FK}{\longleftrightarrow} \partial B(||x-y||/3)) P(0
\stackrel{FK}{\longleftrightarrow} \partial B(||x-y||/3)).
\end{equation}
\end{lemma}

\noindent{\bf Proof of Lemma~\ref{lemma-upper-bound}.}
\begin{eqnarray} \nonumber
\tau_c (y-x) & \leq & P(x \stackrel{FK}{\longleftrightarrow} \partial B_x(||x-
y||/3)
\text{ and } y \stackrel{FK}{\longleftrightarrow} \partial B_y(||x-y||/3)) \\
& \leq & P(x \stackrel{FK}{\longleftrightarrow} \partial B_x(||x-y||/3)
\, | \, y \stackrel{FK}{\longleftrightarrow} \partial B_y(||x-y||/3))
\nonumber \\
&      & P(y \stackrel{FK}{\longleftrightarrow} \partial B_y(||x-y||/3))
\nonumber \\
& \leq & P_{\partial B(||x-y||/3)}^W (0 \stackrel{FK}{\longleftrightarrow}
\partial B(||x-y||/3))
P(0 \stackrel{FK}{\longleftrightarrow} \partial B(||x-y||/3)), \nonumber
\end{eqnarray}
where in the last step we have used FKG. \fbox{} \\

The next lemma uses RSW bounds~\cite{russo,sewe}, namely, that the
probability
$p_{FK}^a (x;r_1,r_2)$
that there is an open FK $a {\mathbb Z}^2$-circuit in an
$(r_1,r_2)$-annulus centered at $x$
is bounded away from zero and one as $a \to 0$
by constants that depend only on $r_1/r_2$.
In fact, we only need
a lower bound; i.e.,
\begin{equation}
\label{eq-rsw}
\text{for any } x \in {\mathbb R}^2 \text{ and some } 0<r_1<r_2<\infty, \,
\liminf_{a \to 0} p_{FK}^a (x;r_1,r_2) \, > \, 0 \, .
\end{equation}
%and for any single choice of $r_1 < r_2$.)
This is not immediate in the Ising case, since there is not
currently a direct proof of RSW for critical FK percolation (as
opposed to the independent percolation case). However, as we explain
after the proof of the lemma, RSW follows from announced results about the
scaling limit of spin cluster boundaries~\cite{smirnov-ising4},
combined with the Brownian loop soup representation of
CLE$_3$~\cite{lw,shw,werner3}; also, the lower bound~(\ref{eq-rsw}) for
some $r_1,r_2$ implies both upper and lower bounds for all $r_1,r_2$.
\begin{lemma} \label{lemma-lower-bound}
%For all percolation measures satisfying FKG and RSW, there
Assuming~(\ref{eq-rsw}), there exists a constant $K>0$ such that
\begin{equation} \nonumber
\tau_c(y-x) \geq K P_{\partial B(||x-y||/3)}^W
(0 \stackrel{FK}{\longleftrightarrow} \partial B(||x-y||/3))^2.
\end{equation}
\end{lemma}

Before giving the proof, we state an immediate consequence of this lemma, the
preceding one and the fact that
$P_{\partial B(r)}^W (0 \stackrel{FK}{\longleftrightarrow} \partial B(r))
\geq P(0 \stackrel{FK}{\longleftrightarrow} \partial B(r))$.
\begin{corollary} \label{cor-comparable}
%For all percolation measures satisfying FKG and RSW,
Assuming~(\ref{eq-rsw}),
$P(0 \stackrel{FK}{\longleftrightarrow} \partial B(||x-y||/3))$ and
$P_{\partial B(||x-y||/3)}^W (0 \stackrel{FK}{\longleftrightarrow} \partial B
(||x-y||/3))$
are comparable (up to constants) as $||x-y|| \to \infty$.
\end{corollary}
%\bigskip
\noindent {\bf Proof of Lemma~\ref{lemma-lower-bound}.}
Let $A(x,r)$ denote the intersection of ${\mathbb Z}^2$
and the annulus with outer radius $r$ and inner radius $r/2$ centered at $x$
and let $circ_{FK}(A(x,r))$ denote the
event that there is an open FK circuit in $B_x(r)$ surrounding
$B_x(r/2)$. Let $F(x,r)$ be the event that $circ_{FK}(A(x,r))$ occurs and the
outermost open FK circuit contained in $B_x(r)$ and surrounding $B_x(r/2)$
is connected to $x$ by an open FK path. We have
\begin{eqnarray} \nonumber
\tau_c (y-x) & \geq & P(F(x,||x-y||/3) \cap F(y,||x-y||/3)
\\ & & \hskip0.5cm \cap \partial B_x(||x-y||/6)
\stackrel{FK}{\longleftrightarrow} \partial B_y(||x-y||/6))  \nonumber \\
          & \geq & K'' P(F(0,||x-y||/3))^2, \label{eq-bound}
\end{eqnarray}
where the second inequality follows from FKG and the constant $K''>0$ follows
from RSW.
%(which remains to be shown).
We then note that
\begin{eqnarray} \nonumber
P(F(0,r)) & = & \sum_{\text{circuits } \gamma} P(\gamma \text{ is the
outermost open
circuit in  $A(0,r)$ and } 0 \stackrel{FK}{\longleftrightarrow} \gamma) \\
          & = & \sum_{\text{circuits } \gamma}P(0 \stackrel{FK}
{\longleftrightarrow} \gamma \, | \,
          \gamma \text{ is the outermost open circuit in } A(0,r)) \nonumber \\
          & & \hskip1.3cm P(\gamma \text{ is the outermost open circuit in } A
(0,r)) \nonumber \\
          & \geq & \sum_{\text{circuits } \gamma}
      P_{\partial B(r)}^W (0 \stackrel{FK}{\longleftrightarrow} \gamma)
          P(\gamma \text{ is the outermost open circuit in } A(0,r)) \nonumber
\\
          & \geq & P_{\partial B(r)}^W (0 \stackrel{FK}{\longleftrightarrow}
\partial B(r)) \nonumber \\
          & & \sum_{\text{circuits } \gamma}
      P(\gamma \text{ is the outermost open circuit in } A(0,r)) \nonumber \\
          & = & P_{\partial B(r)}^W (0 \stackrel{FK}{\longleftrightarrow}
\partial B(r))
          P(\exists \text{ an open circuit in } A(0,r)). \nonumber
\end{eqnarray}
Inserting this bound into~(\ref{eq-bound}) concludes the proof.
%(except for explaining the validity of RSW here).
\fbox{} \\

\noindent{\bf RSW for FK percolation on ${\mathbb Z}^2$.}
In order to get RSW, we assume (from~\cite{smirnov-ising4}) that the
``full" scaling limit of the Ising model converges to (the nested
version of) CLE$_3$. We can then use the representation of CLE$_3$
in terms of the Brownian loop soup~\cite{lw,shw,werner3}, assuming
that $\kappa=3$ corresponds to a density of the Brownian loop soup
below its critical density (which should correspond to $\kappa=4$).
A single Brownian loop has positive probability of ``surrounding" a
disc of fixed radius $r_1$ centered at the origin. Let $\gamma$ be
such a loop and consider its loop-cluster, built recursively from
the (countably many) Brownian loops by saying that any two loops
which touch (and thus cross, with probability one) are in the same
loop-cluster. Given that the density of the Brownian loop soup is
assumed to be below the critical density, the loop-cluster of
$\gamma$
%formed by touching sequences of the (countably many) Brownian loops is
%built recursively by adding to it all the (countably many)
%Brownian loops that touch it,
is contained with probability one inside a sufficiently large disc.
%This means that, for a
%suitably chosen annulus,
Thus, for some $r_2 > r_1$, there is strictly positive
probability that the $(r_1,r_2)$-annulus centered at the origin
contains a CLE$_3$ circuit. Back on the lattice $a {\mathbb Z}^2$,
this gives a positive probability, bounded away from zero as $a \to
0$, that the external boundary of an Ising spin cluster provides
such a circuit. But this in turn implies the same for a closed
(dual) FK circuit, and hence by self-duality at the critical point,
the same for an
%this gives a probability bounded away from both zero and one for an
open FK circuit. %--- i.e., it yields the RSW property.

To conclude the
discussion of RSW, we note that it is not difficult to show
that one can use FKG to obtain from open circuits contained in overlapping
$(r_1,r_2)$-annuli a ``necklace'' structure that provides open crossings
of rectangles of arbitrary aspect ratio (see~\cite{BC08} for more details
about such arguments).
The rectangle crossings can then be used, once again with the
help of FKG, to obtain circuits inside arbitrary annuli with
probability bounded away from zero.
By self-duality one also has closed (dual) crossings of rectangles and
these can be used to bound the probability of open circuits away from one.
%\fbox{} \\

\begin{proposition}
\label{lemma-crossprob} For $z \in {\mathbb R}^2$, let
$N^a(z,r_1,r_2)$ denote the number of distinct clusters $C_i^a$ that
include sites in both $\{y \in a {\mathbb Z}^2: ||y-z|| < r_1\}$ and
$\{y \in a {\mathbb Z}^2: ||y-z|| > r_2\}$.
Assuming~(\ref{eq-rsw}),
for any
$0<r_1<r_2<\infty$, there exists $\lambda \in (0,1)$ such that for
all $z \in {\mathbb R}^2$ and all small $a>0$ and any $k=1,2,\dots$,
\begin{equation}
\label{eq:crossprob}
P(N^a(z,r_1,r_2) \geq k) \, \leq \, \lambda^k .
\end{equation}
It follows that for any bounded $\Lambda \subset {\mathbb R}^2$ and
$\varepsilon>0$, the number of distinct clusters $C_i^a$ of diameter
$> \varepsilon$ touching $\Lambda$ is bounded in probability as $a
\to 0$.
\end{proposition}
\noindent {\bf Proof of Proposition~\ref{lemma-crossprob}.} The proof is
by induction on $k$. For $k=1$, the result follows from RSW since
$N^a(z,r_1,r_2) \geq 1$ is equivalent to the {\it absence\/} of a
closed (dual) circuit in the $(r_1,r_2)$-annulus about $z$,
which by self-duality at the critical point has the same probability
as absence of an open circuit, which in turn is bounded away from
one as $a \to 0$. Now suppose $N^a(z,r_1,r_2) \geq k-1$. Then one
may do an exploration of the $C_i^a$'s that touch $\{y \in a
{\mathbb Z}^2: ||y-z|| < r_1\}$ until $k-1$ are found that reach
$\{y \in a {\mathbb Z}^2: ||y-z|| > r_2\}$, making sure that all
cluster explorations have been fully completed without
obtaining information about the outside of the clusters.
%within the
%$(r_1,r_2)$-annulus.
At that point, the complement $D$ of some
random finite $D^c \subset a{\mathbb Z}^2$ remains to be explored
and the (conditional) FK distribution in $D$ is $P_{\partial D}^F$
with a {\it free\/} boundary condition on the boundary (or
boundaries) between $D$ and $D^c$. By RSW, the $P_{\partial
D}^F$-probability of an open crossing in $D$ of the
$(r_1,r_2)$-annulus is bounded above by
the original $P(N^a(z,r_1,r_2) \geq 1)$. Thus we have
\begin{eqnarray} \nonumber
P(N^a(z,r_1,r_2) \geq k) & = & P(N^a(z,r_1,r_2) \geq k-1) \\ \nonumber
&   & P(N^a(z,r_1,r_2) \geq k|\,N^a(z,r_1,r_2) \geq k-1) \\
& = & P(N^a(z,r_1,r_2) \geq k-1) \, E(P_{\partial D}^F(N^a(z,r_1,r_2) \geq 1))
\nonumber \\
& \leq & P(N^a(z,r_1,r_2) \geq k-1) \, P(N^a(z,r_1,r_2) \geq 1)
\nonumber \\
%& \leq & (P(N^a(z,r_1,r_2) \geq 1))^k \, \leq \, \lambda ^k \, .
& \leq &  \lambda ^k \, . \nonumber
%\nonumber \\
\end{eqnarray}
The last claim of the proposition follows from~(\ref{eq:crossprob}) because
one may choose $O([diam(\Lambda)/\varepsilon]^2)$ points $z_\ell$ in
${\mathbb R}^2$ so that any $C_i^a$ of diameter $> \varepsilon$
touching $\Lambda$ will be counted in $N^a(z_\ell, \varepsilon/4,
\varepsilon/2)$
for at least one $z_\ell$.
\fbox{}
\\

We  next explain how to verify Hypothesis~\ref{hypo} for critical
Ising-FK percolation and for independent percolation; we do not
have a verification for Ising spin percolation, although we expect
it to be true in that case also. It may be of interest to note that
for critical independent percolation,
%(e.g., on ${\mathbb T}$), one
one can obtain a representation like~(\ref{eq:eucfield}), but with
the SLE$_{16/3}$-based measures $\mu_j^{FK}$ replaced by
SLE$_6$-based ones $\mu_j^{IN}$, for the scaling limit of the
lattice ``divide and color'' model~\cite{Hagg01}. (Here
and below we use the letters IN to distinguish independent from FK
percolation.) The original divide and color model, and
the one most analogous to the FK representation of the ($h=0$) Ising
model, takes the open clusters of independent {\it bond\/}
percolation, e.g., on ${\mathbb Z}^2$, and colors them with random
$\pm 1$ signs to define the divide and color spin variables.
In Sec.~\ref{sec-discussion}, we will consider this model as the
density $p$ of open bonds approaches its critical value (we note
that in~\cite{BCM07} a different phase transition is studied). For
independent {\it site\/} percolation, e.g., on ${\mathbb T}$, with
say probability $p$ and $1-p$ for white and black sites, the option
for defining the divide and color spin variables that we will use is
to ``color'' {\it both\/} the white and black clusters with random
signs. It is unclear whether the scaling limit of the critical
($p=1/2$ on ${\mathbb T}$) divide and color model corresponds to
some known conformal field theory. Note that the limit, in terms of
boundaries between clusters of different colors, is a sort of
``dilute CLE$_6$'' and is conformally invariant, but is not itself
described by a CLE since the divide and color model lacks the
``domain Markov property.''
\\

\noindent{\bf Hypothesis~\ref{hypo} for FK percolation on ${\mathbb
Z}^2$.} In this case, the behavior of the two-point function is
known exactly along the $(1,1)$ direction from Ising
calculations, which yield $\tau_c(y-x) \sim K ||x-y||^{-1/4}$
(e.g., \cite{wu}, referred to in~\cite{tracy},
and Chap.~XI of~\cite{mw-book}).
Using Lemmas \ref{lemma-upper-bound} and
\ref{lemma-lower-bound} one obtains that, up to constants, the
two-point function has the same behavior in all directions.
Hypothesis~\ref{hypo} is then satisfied with $2\theta=1/4$.
\\

\noindent{\bf Hypothesis~\ref{hypo} for independent percolation.} From
the analogues of Lemmas \ref{lemma-upper-bound} and
\ref{lemma-lower-bound} for independent percolation, we know that
$\tau_c^{IN}(z)$ is comparable (up to constants) with
$[P(0 \stackrel{IN}{\longleftrightarrow} \partial B(||z||/3))]^2$.
This immediately gives the desired upper bound for $\tau_c^{IN}(x)$. For
the lower bound, it
suffices to show that $P(0 \stackrel{IN}{\longleftrightarrow}
\partial B(r)) \geq K''' (\varepsilon')^{\theta} P(0
\stackrel{IN}{\longleftrightarrow}
\partial B(\varepsilon' r))$ for some
constant $K'''>0$. Using FKG, we have
\begin{eqnarray} \nonumber
P(0 \stackrel{IN}{\longleftrightarrow} \partial B(r)) & \geq &
P(0 \stackrel{IN}{\longleftrightarrow} \partial B(\varepsilon' r) \cap
\partial B(\varepsilon' r/2) \stackrel{IN}{\longleftrightarrow}
\partial B(r) \cap circ_{IN} (A(0,\varepsilon' r))) \nonumber \\
  & \geq & P(0 \stackrel{IN}{\longleftrightarrow} \partial B(\varepsilon' r))
  P(\partial B(\varepsilon' r/2) \stackrel{IN}{\longleftrightarrow} \partial B(r))
  P(circ_{IN} (A(0,\varepsilon' r))) \nonumber \\
  & \geq & K''' (\varepsilon')^{\theta}
  P(0 \stackrel{IN}{\longleftrightarrow} \partial B(\varepsilon' r)), \nonumber
\end{eqnarray}
where $P(circ_{IN} (A(0,\varepsilon' r)) \geq \tilde K$ from RSW, and
$P(\partial B(\varepsilon' r/2) \stackrel{IN}{\longleftrightarrow} \partial B
(r))
\geq \hat K(\varepsilon')^{\theta}$
with $\theta = \alpha-\delta$ for any $\delta>0$, and $\alpha$ denoting
the one-arm exponent.

For site percolation on the triangular lattice, $\alpha$
has been proved to be $5/48$ using SLE computations~\cite{lsw5}. For
other percolation models (e.g., bond percolation on the square
lattice), the five-arm exponent is known to be equal to $2$ (see
Lemma~2 of~\cite{ksz}, Corollary~A.8 of~\cite{ss} and Section~5.2
of~\cite{nolin}), since it can be derived via a general argument that
does not use SLE. Using this and 
the BK inequality~\cite{BK},
we obtain an upper bound of $2/5$ for $\alpha$.
We note, as pointed out to us by P.~Nolin, that a more elementary
argument from~\cite{BK} is available showing that $\alpha \leq 1/2$
without use of the five-arm exponent --- see Eqn.~(2.5) 
of~\cite{chayesnolin}.

\section{Discussion} \label{sec-discussion}

In this paper, we provided a representation
(see~(\ref{eq:eucfield})) for the scaling limit Euclidean random
field $\Phi^0$ associated with the $d=2$ Ising model at its critical
point ($T=T_c$, $h=0$). This field, one of the basic objects of
conformal field theory, is the scaling limit of the magnetization
field $\Phi^a$ (see~(\ref{eq-discrete-field})) on $a {\mathbb Z}^2$
as $a \to 0$. $\Phi^0$ is represented as a sum $\sum_j \eta_j
\mu_j^{FK}$ with random signs $\eta_j$ and finite measures
$\mu_j^{FK}$ that are the limits of rescaled area measures
associated with the macroscopic Ising-FK clusters. These measures
are supported on continuum clusters whose outer boundaries are
described by CLE$_{\kappa}$ loops with $\kappa = 16/3$. A key to
the representation is that natural field strength rescaling
(see~(\ref{eq1-Theta}), (\ref{eq2-Theta})) insures that for bounded
$\Lambda \subset {\mathbb R}^2$, $\sum_j [\mu_j^{FK} (\Lambda)]^2 <
\, \infty$ and hence $\sum_j \eta_j \mu_j^{FK}(\Lambda)$ is convergent
(in $L_2$).

We explained, at the end of
Section~\ref{sec-ising-field-theory}, why the limits $\mu_j^{SP}$ of
area measures for Ising spin clusters do not appear useful for
representing $\Phi^0$. We also noted, towards the end of
Section~\ref{sec-area}, that a field $\sum_j \eta_j \mu_j^{IN}$ can
be constructed using critical clusters from independent in place of
FK percolation, but that its physical significance is unclear.
%A similar representation for the spin field should presumably be valid
%for the 3- and 4-state Potts models, with the Ising-FK clusters and
%area measures replaces by FK cluster with $q=3$ and 4 and their area
%measures and with the $eta_j$'s replaced by variables taking 3 and 4
%values with equal probability, respectively.
We next discuss how the representation~(\ref{eq:eucfield}) could be
extended to off-critical models.

{\it Independent percolation, $p \neq p_c$.\/}
%Can we construct scaling limit area measures for the near-critical
%modification of $CLE_6$ obtained (for site percolation on ${\mathbb T}$)
%by the approach of \cite{cfn1,cfn2} and Garban, Pete and Schramm
%(see~\cite{garban}). Here $p \to p_c$ appropriately as $a \to 0$.
%This should be doable.
Here the percolation density $p$ (say of the white sites on
${\mathbb T}$) converges to the critical density, $p_c$ ($=1/2$),
appropriately as $a \to 0$. In this case, the representation of the
near-critical field should involve area measures from the
near-critical modification of CLE$_6$ obtained by the approach
of~\cite{cfn1,cfn2} and the results of Garban, Pete and Schramm~\cite{gps,garban}.
A feature of that work, which
%we expect to
might also be valid
in the FK-Ising context, is a natural probabilistic
coupling, based on a Poissonian marking of certain
pivotal locations, so that the one-parameter family
of near-critical models parametrized by the strength of the
off-critical perturbation lives on a single probability space.
The appropriate speed at which $p \to p_c$ is
such that the correlation length is bounded away from zero and
infinity. Since it is proved in~\cite{kesten} (and~\cite{nolin}
--- see Theorem~$26$ there) that
%Kesten's scaling relations~\cite{kesten}
%are satisfied and thus
crossing probabilities and multi-arm probabilities are comparable
(up to constants and up to distances of the order of the correlation
length) to those of the critical system,
%(see, e.g., Theorem~26
%of~\cite{nolin}). As a consequence,
RSW still holds and both Hypothesis~\ref{hypo}
and Prop.~\ref{lemma-crossprob} can be verified.
One of the area measures $\mu_j^{IN}$,
corresponding to the unique infinite cluster
(white for $p \downarrow p_c$ and black for $p \uparrow p_c$),
will now have unbounded support and infinite mass over
${\mathbb R}^2$, but its mass in any bounded region $\Lambda$
will be finite and $\sum_j \eta_j \mu_j^{IN}(\Lambda)$ will be
convergent (in $L_2$, at least for a bounded $\Lambda$ not chosen
in a way that depends on $\{\eta_j\}$). We note that for {\it bond\/}
percolation on ${\mathbb Z}^2$ with $p$ the density of open bonds,
it will only be in the $p \downarrow p_c$ near-critical model that
some $\mu_j^{IN}$ has infinite mass.

{\it FK percolation, $T \neq T_c$.\/} Here one keeps $h=0$ in the
Ising model but lets $T \to T_c$ appropriately as $a \to 0$, which
is the analogue of $p \to p_c$ in independent percolation. It is
natural to expect that the claims made above for independent
percolation still hold in this case. Since one is now considering
\emph{bond} (FK) percolation, only for $T \uparrow T_c$ (the analogue
of $p \downarrow p_c$ in independent percolation) will there be an
infinite mass $\mu_j^{FK}$ with unbounded support in ${\mathbb R}^2$.
Including a random sign $\eta_j$ for the infinite mass $\mu_j^{FK}$
means one is taking the scaling limit of the symmetric mixture of
the plus and minus Gibbs measures.

%The most natural approach to this type of off-critical FK percolation
%would be to mimic that for independent percolation of \cite{cfn1,cfn2}
%and Garban, Pete and Schramm (see~\cite{garban}) obtaining pivotal
%(bond rather than site for FK) measures at criticality in the scaling
%limit. The next item would then be to mimic the situation for
%independent percolation, to obtain the corresponding Ising model
%off-critical Euclidean field.

{\it FK percolation, $h \neq 0$.\/} Here one sets $T=T_c$ with $h
\neq 0$ and then lets $h \to 0$ appropriately as $a \to 0$.
Intuitively, in the scaling limit, this should involve formally
multiplying the measure describing the critical continuum system
by a factor proportional to $\exp{(\lambda \int_{{\mathbb R}^2}
\Phi^0(z)\,dz)}$.
%Of the many questions one can then formulate,
%we give one that should have a positive answer.
According to~(\ref{eq:eucfield}) the critical Euclidean field is
given by the sum of all the $\eta_j \mu_j^{FK}$.
Let $\nu^{FK}$ denote the marginal distribution of the process
$\{\eta_j \mu_j^{FK} \}$ of finite signed measures in the plane,
$1_L$ denote the indicator function of the $L \times L$ square
in~${\mathbb R}^2$, $Z_L = \int \exp{(\lambda \Phi^0(1_L))}
d\nu^{FK}$, and finally $d\nu_L^\lambda = (Z_L)^{-1}
\exp{(\lambda \Phi^0(1_L))} d\nu^{FK}$.
We ask: does $\nu_L^\lambda$ converge to some $\nu^\lambda$ as
$L \to \infty$ and is $\Phi^\lambda$, obtained from $\nu^\lambda$
as the sum of its individual signed measures, the physically
correct near-critical Euclidean field?
%Presumably, yes.
Heuristically, the correct normalization to obtain a nontrivial near-critical
scaling limit is such that the correlation length $\xi$ remains bounded away
from zero and infinity. Since $\xi \sim h^{-8/15}$ for small $h$,
this gives $h \sim a^{15/8}$,
which coincides with the normalization needed to obtain
a nontrivial Euclidean field, as can be seen from~(\ref{eq1-Theta})
and the asymptotic behavior of~$\tau_c$.
Using this observation and the
$d=2$ Ising critical exponent $\delta=15$
for the magnetization (i.e., $M \sim h^{1/15}$), the rough computation
(where $\sum_x^L$ denotes the sum over $x$ in $\Lambda_{L /a}$),
\begin{equation} \nonumber
\frac{\langle a^{15/8} \sum_x^1 S_x \exp(a^{15/8} \sum_x^L S_x) \rangle_c}
{\langle \exp(a^{15/8} \sum_x^L S_x) \rangle_c}
\sim a^{-1/8} M(h=a^{15/8}) \sim a^{-1/8} (a^{15/8})^{1/15} = 1,
\end{equation}
suggests a positive answer to the previous questions.

We conclude this section with brief discussions of the
applicability of our approach to higher dimensions, $d>2$,
and to $q$-state Potts models with $q>2$. Although the $d=2$
scaling limit Ising magnetization field $\Phi^0$ should be conformal
with close connections to CLE$_{16/3}$, as we have indicated, very
little conformal or SLE machinery was actually used in our analysis.
Basically, the two main ingredients were (see
Hypothesis~\ref{hypo}) that $\tau_c(y-x)$ behaves at long distance
like $||y-x||^{-\psi}$ with $\psi<d$ and (see
Prop.~\ref{lemma-crossprob}) that as $a=1/{L'} \to0$,
\begin{equation} \label{eq-crossprob2}
P(N^a(0,r_1,r_2) \geq 1) \, = \, P(B(r_1 L')
\stackrel{FK} {\longleftrightarrow} B(r_2 L')^c) \,
\leq \, \lambda \, < \, 1 \, .
\end{equation}
Although such decay of $\tau_c$ should be valid for all $d \geq 2$,
the crossing probability bound~(\ref{eq-crossprob2}) is a different
matter and presumably fails above the upper critical dimension
(see Appendix~A of~\cite{aizenman1}).
%for a related discussion in an independent percolation context).
When it fails, there can be
infinitely many FK clusters with diameter greater than $\varepsilon$
in a bounded region and so Prop.~\ref{lemma-zero-limit} would not
preclude $\Phi^0$ from being a Gaussian (free) field. But it
appears that at least for $d=3$, both~(\ref{eq-crossprob2}) and
a representation of $\Phi^0$ as a sum of finite measures with
random signs ought to be valid.

As pointed out to us by J.~Cardy, an analogous representation for the
scaling limit magnetization fields of $q$-state Potts models also ought
to be valid, at least for values of $q$ such that for a given $d$, the
phase transition at $T_c$ is second order. The phase transition is
believed to be first order for integer $q \geq 3$ when $d \geq 3$
and for $q > 4$ when $d = 2$ --- see~\cite{wu-review}; this leaves,
besides the Ising case, $d = 2$ and $q = 3$ and~$4$. We denote the
states or colors by $1,2,\dots,q$ and recall that in the FK representation
on the lattice, all sites in an FK cluster have the same color while the
different clusters are colored independently with each color equally
likely. In the scaling limit, there would be finite measures
$\{\mu_j^{FK,q}\}$ and the magnetization field in
the color-$k$ direction would be $\sum_j \eta_j^k \mu_j^{FK,q}$
with the $\eta_j^k$'s taking the value $+1$ with probability $1/q$
(for the color $k$) and the value $-1/(q-1)$ with probability $(q-1)/q$
(for any other color). For a fixed $k$ the $\eta_j^k$'s would be
independent as $j$ varies, but for a fixed $j$ they would
be {\it dependent\/} as $k$ varies because $\sum_k \eta_j^k =0$.

\bigskip

{\bf Acknowledgements.} The authors thank the Centre
de Recherches Math\'ematiques, Montr\'eal, for hospitality
during August, 2008
%``Stochastic Loewner Evolution and Scaling Limits,"
and the Institut Henri Poincar\'{e} - Centre Emil Borel,
as well as Univ.~Paris Sud [C.M.N.] and
\'{E}cole Normale Sup\'{e}rieure [F.C.],
for hospitality in Paris during
October and November, 2008.
%part of the 2008-2009 Thematic Program on ``Probabilistic Methods in
%Mathematical Physics."
C.M.N. thanks the
department of mathematics of the Vrije Universiteit Amsterdam
for its hospitality during a visit in 2007, when
the present work was started and during a visit in 2008. F.C. thanks
the Courant Institute of Mathematical Sciences for its hospitality
during two visits in 2008.
The authors thank Douglas Abraham for communications concerning
critical Ising two-point functions, Michael Aizenman, Vincent Beffara,
John Cardy, Oscar Lanford, Pierre Nolin, 
Oded Schramm, Stas Smirnov and Alan Sokal
for useful conversations, Christophe Garban for discussions about his
work with Pete and Schramm, and Wouter Kager for providing Figure~\ref{loopsfigure}.
%the figure we use for loops in the medial lattice.

\end{document}